\documentclass[journal]{IEEEtran}
\usepackage{cite}
\ifCLASSINFOpdf
\else
\fi
\usepackage{graphicx}
\usepackage[reqno]{amsmath}
\usepackage{bm}
\begin{document}
\title{Grid Optimal Integration of Power Ships}
\author{Motab~Almousa,~\IEEEmembership{Student Member,~IEEE}

\thanks{Motab Almousa is with the Department of Electrical and Computer Engineering, Louisiana State University, Baton Rouge, LA, 70803 USA e-mail: (malmou2@lsu.edu)}
}
\markboth{ ~24 May~2019}%
{Shell \MakeLowercase{\textit{et al.}}: Bare Demo of IEEEtran.cls for IEEE Journals}
\maketitle
\begin{abstract}
Power-generating ships or Power Ships (PSs) are considered one of the largest mobile energy resources. In this paper, a model is proposed to evaluate the integration of PSs into power grid operations. The model optimally coordinates the ships to enhance the grid objective; its solution provides optimal generation resource scheduling, as well as optimal scheduling and routing of the ships. IEEE 6-bus and IEEE 118-bus case studies were considered to model the system and to validate the effectiveness of the proposed model.
\end{abstract}
\IEEEpeerreviewmaketitle
\section*{Nomenclature}
\addcontentsline{toc}{section}{Nomenclature}
\subsection*{Sets}
\begin{IEEEdescription}[\IEEEusemathlabelsep
\IEEEsetlabelwidth{$V_1,V_2,V_3$}]
\item[$\mathcal{B}$] Set of buses.
\item[$\mathcal{G}$] Set of generation units.
\item[$\mathcal{L}$] Set of transmission lines.
\item[$\mathcal{P}$] Set of ports.
\item[$\mathcal{S}$] Set of power ships.
\item[$\mathcal{T}$] Set of time periods.
\end{IEEEdescription}

\subsection*{Indices}
\begin{IEEEdescription}[\IEEEusemathlabelsep\IEEEsetlabelwidth{$V_1,V_2,V_3$}]
\item[$i,j$] 		Indices of ports.
\item[$n,m$] 		Indices of buses.
\item[$t,t1,t2$] 	Indices of time.
\end{IEEEdescription}

\subsection*{Parameters}
\begin{IEEEdescription}[\IEEEusemathlabelsep\IEEEsetlabelwidth{$V_1,V_2,V_3$}]
\item[$C_{i,s,t}^{D}$]   	Departure cost for ship $s$ from port $i$ at time $t$.
\item[$C_{i,s,t}^{E}$]    	Entering cost for ship $s$ to port $i$ at time $t$.
\item[$C_{i,s,t}^{W}$]   	Waiting cost for ship $s$ in port $i$ at time $t$.
\item[$C_{i,s,t}^{S}$]   	Sailing cost for ship $s$ in port $i$ at time $t$.
\item[$D_{m,t}$]  			Load Demand of bus $m$ at time $t$.
\item[$F(P)^{g}_{a,b,c}$]   Production cost function of generator $g$.
\item[$F(PS)^{s}_{a,b,c}$]  Generation cost  function of ship $s$.
\item[$M_{g}^{On}$]    		Generator $g$ minimum on time.
\item[$M_{g}^{Off}$] 		Generator $g$ minimum off time. 
\item[$T_{S}^{ij}$] 		Travel time of ship $s$ from port $i$ to port $j$.
\item[$RU_{g}$] 			Maximum ramp-up rate of generator $g$.
\item[$RD_{g}$] 			Maximum ramp-down rate of generator $g$
\item[$SUC_{g}$] 			Start-up cost of unit $g$.
\item[$SDC_{g}$] 			Shutdown cost of unit $g$.
\item[$F_{l}^{Max}$]  		Maximum flow limit of line $l$.
\item[$X_l$]  				Electrical reactance of line $l$.
\item[$\Psi_{n,t}$]  		Shed-load factor in bus $n$ at time $t$.
\end{IEEEdescription}

\subsection*{Variables}
\begin{IEEEdescription}[\IEEEusemathlabelsep\IEEEsetlabelwidth{$V_1,V_2,V_3$}]
\item[$F_{l,t}$]  			Power flow in line $l$ at time $t$.
\item[$P_{g,t}$] 			Output power of generator $g$ at time $t$.
\item[$PS_{s,i,t}$] 		Output power of ship $s$ in port $i$ at time $t$.
\item[$SHD_{m,t}$] 			Load shedding from bus $m$ at time $t$.
\item[$\theta_{t}$] 		Bus phase angle at time $t$.
\end{IEEEdescription}

\subsection*{Binary Variables}
\begin{IEEEdescription}[\IEEEusemathlabelsep\IEEEsetlabelwidth{$V_1,V_2,V_3$}]
\item[$U_{g,t}$]      		1 if Generator $g$ is operating at time $t$.
\item[$SU_{g,t}$]  			1 if Generator $g$ started up at time $t$.
\item[$SD_{g,t}$]  			1 if Generator $g$ is shut down at time $t$.
\item[$V_{i,s,t}$]     		1 if ship $s$ is located in port $i$ at time $t$.
\item[$W_{i,s,t}$]    		1 if ship $s$ is waiting in port  $i$ at time $t$.
\item[$O_{i,s,t}$]     		1 if ship $s$ is operating in port $i$ at time $t$.
\item[$V_{ij,s,t}^{S}$]     1 if ship $s$ is sailing from port $i$ to $j$ at time $t$.
\item[$V_{i,s,t}^{D}$]     	1 if ship $s$ departed from port $i$ at time $t$.
\item[$V_{i,s,t}^{E}$]      1 if ship $s$ entered port $i$ at time $t$.
\end{IEEEdescription}

\newtheorem{Theorem}{Theorem}[section]
\newtheorem{Lemma}{Lemma}[section]
\newtheorem{Proof}{Proof}[section]
\section{Introduction}
\IEEEPARstart{E}{l}electricity plays a major role in our daily lives. Countries with many islands face the challenge of supplying electricity offshore, and ships have become useful for the supply of electricity to such countries. Indonesia alone has 13,000 islands, and supplying them with electricity is problematic due to the high cost of transporting the fuel. Hence, the Japanese and Indonesian governments are considering deploying power-generating ships (PSs) to the islands as a solution \cite{MM}.

Nowadays, PSs are attracting considerable attention as a long-term and short-term energy resources. They have become beneficial in the supply of electricity in many countries. The use of PSs has been proven to be economical and effective in helping many islands and other countries surrounded by bodies of water to have a reliable source of electric power \cite{camacho_2018}. The use of such ships as an energy resource is not a new concept; however, the increase in demand for electricity in combined with advances the manufacturing and technology and economies of scale have led to advances have led increased feasibility.

PSs can be conventional, hybrid, or electrical\cite{21}. Conventional PSs use oil or natural gas as fuel, while hybrid and electrical models have many similarities to other electric vehicles. A power ship's size and speed are essential factors in its deployment\cite{fair}. Large ships are considered static assets, while small ships are considered more appropriate for tactical and operational use. Depending on the prospective task, large ships can be deployed for months, while smaller ships can be used for a day or even a few hours. Nevertheless, wide range of ships and vessels can produce power on the order of 10 MW. That means, they technically can be deployed as power-generating ships if needed\cite{ss}.

\begin{figure}[!t]
\centering
\includegraphics[width=3.0in]{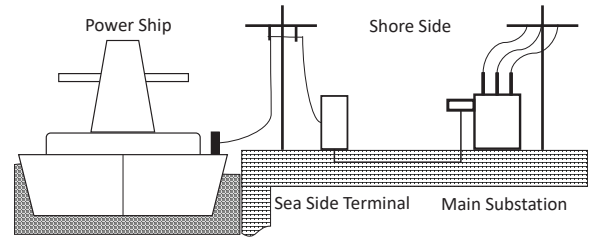} 
\caption{Ship to grid (S2G) system configuration.}
\label{Ship}
\end{figure}

Considering the existing infrastructure, no significant changes are needed to enable the ship-to-grid (S2G) technology. Ports are equipped with high power electrical systems to supply the ships' auxiliary system with electricity needed while moored. An example of a typical schematic S2G configuration is shown in Fig. \ref{Ship}. The major components of the onshore electrical system include high-power cables, switchgear, and transformers, which are naturally capable of working in both directions. This capabily offers a potential for wider ships deployment and integration with the electric system to send energy back to the grid \cite{ss}. 
  
That said, PSs can enhance the power system security, resilience, and economics by providing energy to areas where it's most needed in a fast and reliable manner, e.g. areas affected by natural disasters or grid outages. Moreover, ships with energy storage capabilities and solar panels can play a significant environmental role by catching excess renewable energy and transmitting it to the electric grid of islands. 

A significant focus in the literature is on the ships' electrical system known as the ship power system (SPS). SPS studies include bus reconfiguration, protection, and stability and dynamics of ship electrical systems. For example, a model for an optimal load shedding study of SPS was presented by \cite{mzk02}. A contingency-based method for optimal SPS bus reconfiguration was proposed in \cite{10}. In \cite{7}, different models are presented to assess the dynamics of a ship's primary electrical components. 

Moreover, a remarkable amount of interest has been focused on extending study of SPS studies to  microgrid-based studies \cite{DCMT}. For example, renewable energy integration to the ship has been discussed in \cite{t17,kdbd17,AES1,AES2}. In \cite{t17}, a control algorithm was proposed to optimally integrate photovoltaic (PV) panels on the ship. In \cite{kdbd17}, a ship was modeled with different energy resources, such as wind and solar, with the objective of designing a load frequency controller. Moreover, In \cite{AES1} and\cite{AES2}, an energy management model was presented that included operational security constraints and was suitable for integration with an all-electric ship. 

Most studies have focused on the internal electrical system of the ships as a system that is isolated and independent from the main grid. In fact, very few studies have discussed the interaction between the electric grid and ships. The potential of deploying ships as dedicated PSs to provide power to the electrical grid, has even received less research attention.

On the other hand, the potential impact of energy resource mobility on grid operations even in relatively small-scale integration, such as the vehicle-to-grid concept in electric vehicles (EVs), is not negligible. Transportation models have been widely proposed and discussed in the literature, particularly in relation to, EVs and railway systems. The EV literature is rich with models that solve different kinds of problem, e.g. \cite{EV1,EV2,EV3,EV4,EV5}. In general, the solutions of such models give the optimal energy scheduling (Charging/Discharging), as well as a routing solution \cite{awcgsj17}. Different objectives can be considered, such as minimum grid cost, minimum EV owner (or fleet) cost, or minimum average street congestion. 

In work with important implications for this paper, the authors of \cite{kws12} have demonstrated the potential of EV fleets coordination on the operational costs of power system. Their study highlighted the economic feasibility of EVs if their deployment is coordinated with grid operations and takes into consideration the major transportation and electric grid modeling requirements. 

Moreover, in \cite{1,2,2s}, the authors proposed a model to consider railway networks to enhance  grid operations. The primary objective of their model was to utilize the railway network to transport energy between different locations; The results shows that carrying batteries between train stations can relieve grid congestion and network bottlenecks. Their proposed model demonstrates a significant economic savings potential based on real-world assumptions.
 
However, to the best of our knowledge, ships have not been addressed in the literature as a solution to enhance electric grid operations and economics. Unfortunately, models to integrate EVs into the grid cannot be directly applied to maritime-based problems. Ships, ports, and maritime transportation models have fundamental differences from EVs models. Ships' sailing, routing, and ports management are examples of areas that require explicit models for the problem. The development of a model that is practical, efficient, and incorporates power ships into power system operations is therefore crucial.
In this paper, a maritime model is introduced to coordinate and adequately incorporate power-generating ships into grid operations. The benefits of the proposed Maritime-based Energy Scheduling and Coordination (MESC) Model include:
\begin{itemize}
\item[$1.$] {Enabling better estimation of the economic and operational potential of incorporating PS's into power grids.}\\
\item[$2.$] {Providing quantitative and qualitative measures of the impact of the intrinsic characteristics of maritime transportation systems on power system operations.}\\
\item[$3.$] {Addressing some unseen challenges, such as the computational complexity of such a model.}
\end{itemize} 
The remaining of the paper is organized as follows: The second section describes the problem and the proposed model. Section III outlines the solution methodology. In sections IV and V, the mathematical formulation is provided. In section VI, case studies are presented and solved. A discussion on the computational complexity of the problem is provided in section VII.  Conclusions are drawn in the final section.

\section{Maritime-based Energy Scheduling and Coordination Model}

The proposed model as defined includes an energy generation scheduling problem and a maritime-based transportation problem. The solution should determine the optimal routes and schedules for the PSs, including deciding their hourly power dispatch for optimal grid operation. This plan must satisfy the operational and technical constraints of the ship in addition to the grid operational constraints of the grid. Different maritime-based considerations were adopted from real-world maritime transportation models presented and discussed in \cite{OR1,OR2,OR3,OR4}, including ports management, e.g., operations capacity and berth limits. A discrete time space was considered due to the nature of the problem. Furthermore, a deterministic model was assumed for modeling the problem, including travel time between the ports. Generally, there were no operational restrictions on travel between ports; However, we assumed that PSs would not be allowed to be traveling or sailing at the end of the planning horizon. Despite this restriction, they could be located at any port. Sailing time included all setup time needed to accommodate the ship. Furthermore, a ship's power generation was modeled by the same constraints imposed on thermal units inland, such as ramping rates and minimum on and off times. Moreover, they were required to operate upon the arrival at their destination. A One-time cost upon the departure and arrival was imposed. Additionally, sailing costs were considered depending on the ship type and distance to be traveled. Ships were not required to operate during the planning time horizon; they could stay idle. However, since ships usually wait or moor in ports, a waiting cost was imposed during the idle time. 

On the other hand, generation resources scheduling with grid constraints was considered as the electrical grid operation problem. The term grid-constrained unit commitment (GCUC) was introduced to refer to the problem. Since the proposed energy problem is different than the most commonly known, security constraint unit commitment (SCUC), a new terminology was needed. In the SCUC problem, the grid is operated by an independent system operator (ISO) that coordinates the energy resources while taking into account specific security operational requirements imposed by a higher regulatory authority, such as the Federal Energy Regulatory Commission (FERC) in the United States. However, offshore, grid operations are usually undertaken locally by the grid operator, under the grid operator or local authorities requirements. 

The problem was mathematically formulated as a mixed integer programming (MIP) problem.

\section{MESC Solution Methods and Approaches}

\begin{figure}[!t]
\centering
\includegraphics[width=3in]{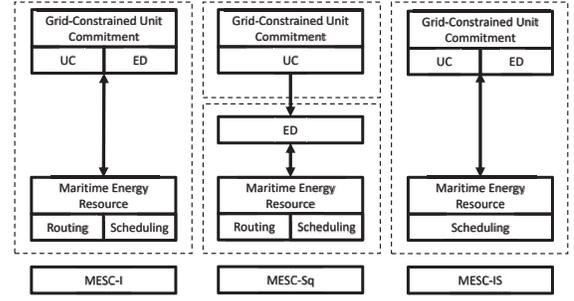} 
\caption{Different solution approaches to solve the problem.}
\label{Approach}
\end{figure}

While the MESC model can be solved for small problems, the problem becomes much harder for larger cases. Therefore, in addition to the original or exact solution approach, where the problem is solved all at once, we considered two approximate solution approaches to the model. We refer to the general solution of the model by the term integrated solution approach, or MESC-I (I referring to integration). This approach is, of course, the optimal way to solve the model. A common solution approach in combinatorial problems is to solve the problem in two sequential steps, usually suitable when the original problem is naturally consisting of two decomposable sub-problems, as in this case. For our model, the solution approach is as follows: First, the energy scheduling problem (GCUC) is solved and the binary variables are fixed, then the complete model is solved again. It is trivial to see a significant computational time reduction, since the number of the binary variables in each step is reduced, even though the problem is solved twice instead of once. Since the problem is solved in a sequence in this approach, we will refer to this solution approach as MESC-Sq, (Sq for sequential). Moreover, another intuitive approach is when the grid operator is only responsible for the ship's operations. Assuming ships routing is part of a different problem with longer time span, for example, the weekly operational planning with time spans of days and ships' travel is restricted to specific days in the week. The term MESC-IS is introduced to refer to this solution approach, (I and S referring to integrated and stationary within the problem time horizon). Fig. \ref{Approach} shows the proposed solution approaches.

\section{MESC Mathematical Formulation I}
The mathematical formulation of ships routing and coordination problem is provided in this section.
\subsection{Power Ships Objective Function}
\begin{equation*}
\begin{split}
Min \sum_{ij \in \mathcal{P}^{+}} \sum_{s \in \mathcal{S}} \sum_{t\in \mathcal{T} } & C_{i,s,t}^{E}  V_{i,s,t}^{E}+ C_{i,s,t}^{D} V_{i,s,t}^{D} +\\
& C_{i,s,t}^{W}  W_{i,s,t}+  C_{ij,s,t}^{S} V_{ij,s,t}^{S} +\\
& F(PS)^{s}_{a,b,c} + SUC_{s}SU_{s,t} +SDC_{s}SD_{s,t}\\
\end{split}
\end{equation*}
The total cost is categorized into 5 parts as follow. $C_{ij,s,t}^{S}$ is the cost of sailing, which includes fuel and other costs incurred during the travel. Ships entering and departing costs are given by $C_{i,s,t}^{D}$ and $C_{i,s,t}^{E}$, respectively. Departure cost $C_{ij,s,t}^{D}$ is imposed by ports operators on the ships;  to account for costs incurred during the ships' departure. Entering or set up cost to accommodate the ships is given by $C_{ij,s,t}^{E}$. $C_{ij,s,t}^{W}$ is waiting cost imposed on the ship while it waits or moors at ports. The remaining part is the ships' energy generation costs, where $F(PS)^{s}_{a,b,c}$ is the PS generator quadratic cost function and $SUC_{s}$ and $SDC_{s}$ are the start-up and shutdown costs, respectively. $\mathcal{P}^{+}$ refers to a subset of all feasible routes between the buses. 

\subsection{Power Ships Flow Constraints}
\begin{equation}
 \label{a}
\begin{split} 
 & \sum_{i \in \mathcal{P}} V_{i,s,t}+ \sum_{ij \in \mathcal{P}^{+}} V_{ij,s,t}^{S}= 1  \quad   \forall s \in \mathcal{S},
  t\in \mathcal{T} 
  \end{split}
\end{equation}  
\begin{equation} 
\label{b}
 W_{i,s,t}  + O_{i,s,t} \geq V_{ij,s,t-1}^{S} - V_{ij,s,t}^{S}      \enspace  \forall s \in \mathcal{S},
  i \in \mathcal{P},t \in \mathcal{T}
\end{equation}
\begin{equation}
\label{c} 
 V_{j,s,t} \geq V_{ij,s,t-1}^{S} - V_{ij,s,t}^{S} \quad  \forall s \in \mathcal{S},
  j \in \mathcal{P},t \in \mathcal{T} 
\end{equation}
Constraint \eqref{a} ensures that the ship at any time is either in a port or sailing between ports. Constraints \eqref{b} implies that a ship can only sail out of a port after it was either operating or waiting at that port. Constraint \eqref{c} ensures that ship must enter its destination.

\subsection{Power Ships Arrival-Departure Logic}
\begin{equation} 
\label{d} 
 V_{i,s,t}^{D} \geq V_{i,s,t-1}- V_{i,s,t}      \enspace  \forall i \in \mathcal{P}, 
 s \in \mathcal{S},t \in \mathcal{T}  
\end{equation}
\begin{equation}
\label{e}  
  V_{i,s,t}^{E} \geq V_{i,s,t} - V_{i,s,t-1}    \enspace \forall i \in \mathcal{P},
   s \in \mathcal{S},t \in \mathcal{T}  
\end{equation}
\begin{equation}
\label{f}  
  V_{i,s,t}^{E} +V_{i,s,t}^{D}  \leq 1  \qquad  \forall i \in \mathcal{P}, s \in \mathcal{S},t \in \mathcal{T}  
\end{equation}

Ships entering and departing logic are given in constraints \eqref{d}-\eqref{e}. Constraint \eqref{f} prevents ships from entering and departing a port at the same time period.

\subsection{Power Ships Operational Constraints}
\begin{equation} 
\label{g} 
 W_{i,s,t} = V_{i,s,t} - O_{i,s,t}    \enspace \forall i \in \mathcal{P},
   s \in \mathcal{S},t \in \mathcal{T}  
\end{equation}
\begin{equation}
\label{h} 
 O_{i,s,t} \leq V_{i,s,t-1}    \enspace \forall i \in \mathcal{P},
   s \in \mathcal{S},t \in \mathcal{T}   
\end{equation}
\begin{equation}
\label{i} 
 O_{i,s,t} \geq V_{i,s,t}^{E}  \enspace \forall i \in \mathcal{P},
   s \in \mathcal{S},t \in \mathcal{T}   
\end{equation}

Where, $V_{i,s,t}$ is 1, whenever ship $s$ is located at port $i$ at a given time period $t$. Constraint \eqref{g} ensures that the ship is either operating or waiting at the port.
Constraint \eqref{h} restricts the ship from operating the generation units while sailing. Thus, a ship can only operate if it is at a port. Constraint \eqref{i} implies that the ship must operate upon its arrival.

\subsection{Power Ships Travel Time Constraints}

\begin{multline}
\label{j}
\sum_{\tau \in \Omega} V_{ij,s,\tau}^{S} \leq T_{s}^{ij}  + (1-V_{s,i,t}^{D})M  \\
  t1,t2  \in \mathcal{T} \mid t2 > t1+ T_{s}^{ij},     t1 \leq  \Omega  \leq  t2
\end{multline}
\begin{multline}
\label{k}
\sum_{\tau \in \Omega} V_{ij,s,\tau}^{S} \geq T_{s}^{ij}  - (1-V_{s,i,t}^{D})M  \\
  t1,t2  \in \mathcal{T} \mid t2 > t1+ T_{s}^{ij},     t1 \leq  \Omega  \leq  t2
\end{multline}
\begin{multline}
\label{l}
V_{s,ij,t}^{S}+V_{s,ij,t-T_{s}^{ij}}^{S} \leq 1 \enspace \forall i,j \in \mathcal{P}^{+},\\
 s \in \mathcal{S},t \in \mathcal{T} \mid t> T_{s}^{ij} 
\end{multline}

Constraints \eqref{j}-\eqref{k} impose the travel time between two ports. $T_{s}^{ij}$ is the time needed for ship $s$ to travel from port $i$ to port $j$. Moreover, to make sure that the ship ends its route in a port by the end of the provided travel time, constraint \eqref{l} is imposed. $M$ is a large positive number.

\subsection{Ports Operational Constraints}
Ports are limited to a certain operating capacity. Moreover, Departure or berth is labor intense process and also is limited by a capacity limit. 
The following constraints ensures the limits are implied.
\begin{equation} 
\label{m}
\sum_{s \in \mathcal{S}}  O_{i,s,t}  \leq POC_{i}     \quad \forall i \in \mathcal{P},
   t \in \mathcal{T}  
\end{equation}
\begin{equation} 
\label{n}
\sum_{t \in \mathcal{T}}  V_{i,s,t}^{D}  \leq PDC_{i}    \quad \forall i \in \mathcal{P},
   s \in \mathcal{S}   
\end{equation}
Constraint \eqref{m} guarantees that the number of ships operating at a port do not exceed the port capacity limit $POC_{i}$.
constraint \eqref{n} restricts the number of ships in the port to not exceed the berth capacity of the port $PDC_{i}$.
\subsection{Power Ships Generation Unit Constraints}
PSs generation unit constraints are assumed to be the same as the thermal generation unit constraints provided in next section.
\subsection{Binary and Non-negativity constraints}
\begin{equation} 
\label{o}
 O_{i,s,t}  \in \{ 0,1 \} \enspace \forall i \in \mathcal{P},
   s \in \mathcal{S},t \in \mathcal{T} 
\end{equation}
\begin{equation}
\label{p} 
 W_{i,s,t}   \in \{ 0,1 \}     \enspace \forall i \in \mathcal{P},
   s \in \mathcal{S},t \in \mathcal{T}  
\end{equation}
\begin{equation}
\label{q} 
 S_{i,s,t}   \in \{ 0,1 \}      \enspace \forall i \in \mathcal{P},
   s \in \mathcal{S},t \in \mathcal{T}  
\end{equation} 
\eqref{o}-\eqref{q} are the binary constraints.
\section{MESC Mathematical Formulation II}
The mathematical formulation of the energy resources scheduling (GCUC) problem is provided in this section, (MESC II).
\subsection{Objective Function of GCUC}
The objective function of the energy scheduling sub-problem is: 
\begin{equation*}
 Min \sum_{g \in \mathcal{G}} \sum_{t \in \mathcal{T}} [ F(P)^{g}_{a,b,c} + SUC_{g}SU_{g,t} +SDC_{g}SD_{g,t}  ]
 \end{equation*}
 Where, $F(P)^{g}_{a,b,c}$ is the quadratic cost function given as follow,  
$a_{g}U_{g,t}+b_{g}P_{g,t}+c_{g}P_{g,t}^{2}$. Start up and shut down costs 
of the units are also considered.

\subsection{Generator Limits Constraints}
The generation units are limited by generation capacity limits.
\begin{equation}
\label{r}
 U_{g,t} P_{g}^{Min}  \leq P_{g,t} \leq U_{g,t} P_{g}^{Max}   \quad \forall g \in \mathcal{G}, t \in \mathcal{T} 
\end{equation}
Constraint \eqref{r} implies the upper and lower generation capacity limits.

\subsection{Ramping Rate Constraints}
Generation units are constraint by ramping rate limits. 
\begin{equation}
\label{u}
P_{g,t}-P_{g,t-1} \leq RU_{g} \qquad \forall g \in \mathcal{G},
 t \in \mathcal{T} 
\end{equation}
\begin{equation}
\label{v}
P_{g,t-1}-P_{g,t} \leq RD_{g} \qquad \forall g \in \mathcal{G},
 t \in \mathcal{T} 
\end{equation}
Constraints \eqref{u}-\eqref{v} enforce the ramp-up and ramp-down generation production limits.

\subsection{Start-up and Shutdown Logic} 
The variables to indicate the start-up and shutdown status of the generation units are described by the following equation and constraint.
\begin{equation}
\label{s}
 SU_{g,t} - SD_{g,t} =   U_{g,t} - U_{g,t-1}   \quad \forall g \in \mathcal{G},
 t \in \mathcal{T} 
\end{equation}
\begin{equation}
\label{t}
 SU_{g,t} +SD_{g,t} \leq 1   \qquad \forall g \in \mathcal{G},
 t \in \mathcal{T} 
\end{equation}
In \eqref{s}, the start-up and shutdown logic  is obtained. Constraint \eqref{t} ensures that the units can not start-up and shutdown at the same time.

\subsection{Minimum Up and Down Time Constraints}
Thermal generation units are usually restricted to specific time period before they are allowed to change their operation status.
\begin{alignat}{2}
\label{x}
T_{g,t}^{on} \geq UT_{g} (U_{g,t}-U_{g,t-1})    \quad \forall g \in \mathcal{G},
 t \in \mathcal{T} 
\end{alignat}
\begin{alignat}{2}
\label{y}
T_{g,t}^{off} \geq DT_{g} (U_{g,t-1}-U_{g,t})	\quad \forall g \in \mathcal{G},
 t \in \mathcal{T} 
\end{alignat}
In\eqref{x}-\eqref{y}, the minimum up and down time for the generation units are ensured.

\subsection{Transmission Network Constraints}
Constraints \eqref{z}-\eqref{z5}  represents the electrical network model and grid operational constraints.
\begin{alignat}{2}
\label{z}
F_{l,t} =  \frac{\theta_{m,t}-\theta_{n,t}}{X_{l}}   \quad \forall l \in \mathcal{L},
 ( m,n )\in \mathcal{L}^{+}, \forall t \in \mathcal{T}
\end{alignat}
\begin{multline}
\label{z1}
\sum_{l  \in \mathcal{L}^{n}}  F_{l}  + \sum_{g \in \mathcal{G}^{n}} P_{g,t} + \sum_{s \in S} PS_{s,i,t} -\\
  D_{m,t}+SHD_{m,t}=0 \quad  \forall i \in \mathcal{P} \mid i=m ,t \in \mathcal{T}
\end{multline}
 \begin{equation}
 \label{z2}
0 \leq SHD_{m,t} \leq \Psi_{m,t} D_{m,t} \qquad \forall m \in \mathcal{B},
 t \in \mathcal{T} 
\end{equation} 
\begin{equation}
\label{z3}
 F_l^{MIN} \leq F_{l} \leq F_l^{MAX}  \qquad \forall l \in \mathcal{L}
\end{equation}
\begin{equation}
\label{z4}
\theta^{MIN}  \leq \theta_{m,t}  \leq \theta^{MAX}  \quad \forall   m \in \mathcal{B},t \in \mathcal{T} 
\end{equation}
\begin{equation}
\label{z5}
\theta_{Ref,t}=0    \qquad   \forall  t \in \mathcal{T} 
\end{equation}
In \eqref{z}, the power flow in lines is calculated. The nodal power balance is enforced by constraint \eqref{z1}. Constraint \eqref{z2}, ensures that load shedding is restricted to a predefined limit by the grid operator, e.g. if $Psi_{m,t}$ is equal to 1, the load can be totally shed from that bus. Constraint \eqref{z4}, imposes the phase angle operational limits. Line flow limits are ensured by constraint \eqref{z3}, where $F_l^{MIN} = -F_l^{MAX}$, and similarly, $\theta^{MIN}$ is equal to -$\theta^{MAX}$. Constraint \eqref{z5} provides the reference angle.
\subsection{Binary and Non-negativity Constraints}
\begin{equation} 
\label{z6}
 U_{g,t}  \in \{ 0,1 \}       \enspace \forall g \in \mathcal{G}, t \in \mathcal{T} 
\end{equation}
\begin{equation}
\label{z7} 
 SU_{g,t}   \in \{ 0,1 \}     \enspace \forall g \in \mathcal{G}, t \in \mathcal{T}   
\end{equation}
\begin{equation}
\label{z8} 
 SD_{g,t}   \in \{ 0,1 \}      \enspace \forall g \in \mathcal{G}, t \in \mathcal{T}  
\end{equation} 
\eqref{z6}-\eqref{z8} are the binary constraints.

\section{Case Studies}

\begin{figure}[!t]
\centering
\includegraphics[width=2.8in]{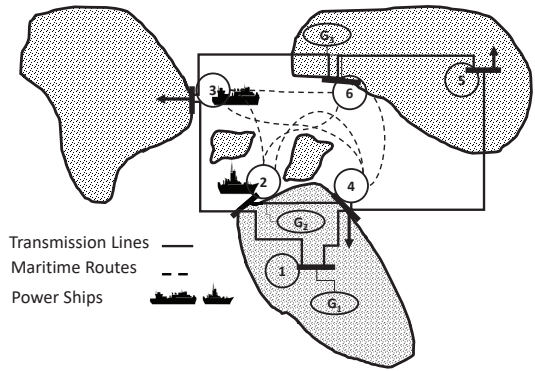} 
\caption{Maritime system with 3 islands and IEEE 6-bus electrical network.}
\label{system}
\end{figure}

The IEEE 6-bus and the IEEE 118-bus systems were utilized to model and solve the problem. To avoid the need to modify the systems, we assumed the three regions or islands, as shown in Fig. \ref{system} for the 6-bus system, are electrically interconnected, yet isolated from the main grid. Ships could access buses 2,3,4, and 6. Buses 2 and 4 were on the same island but at different geographical locations, while buses 3 and 6 were the accessible buses on the other two islands. Case studies II and III were conducted on the 118-bus system for different loading scenarios. However, in both cases, we considered buses 7, 10, 70, 75, 87, and 97 to be the accessible buses or ports. Initially, the ships were on buses 2, and 3 in the 6-bus system, and located at buses 87 and 97 in the 118-bus system. 

As mentioned in the model description, PSs have thermal generation characteristics and limits. Data from the IEEE 118-bus case study were adopted for PSs' generation unit characteristics and operational limits. Two ships of different sizes were considered. The larger ship (PS 1) had the same characteristics as generation unit number 30 in the IEEE 118-bus case, while the smaller ship (PS 2) had the generation characteristics of unit number 42. Table \ref{tgcps}, shows the generation characteristics of the two PSs. 

The larger ship could travel to any of the neighboring ports in 3 hours, while the smaller one needed 2 hours. However, different travel times can be considered without changes to the mathematical formulation of the model. Ports costs are illustrated in table \ref{tpc}. The quadratic term in the generation cost function was ignored to avoid the need for a quadratic solver, but it can be expressed by a linear piece-wise approximation. Without loss of generality, only the linear terms in the cost function were considered. Moreover, the cost for any unmet demanded load at any bus was set to \$1000 per MW. Buses 1 and 10 were selected to be the reference buses for Case I and cases II-III, respectively.

The general algebraic modeling system (GAMS) was used to build the model. The stopping criteria or the termination condition for the solver was set to zero duality gap with a maximum running time of 2500 seconds. In all case studies, the model was run on the NEOS server \cite{NEOS1,NEOS2,NEOS3} with CPLEX as the solver. The NEOS server specifications, include two Intel Xeon X5660 CPUs @ 2.8GHz (12 cores total), and 64GB RAM; More details can be found in \cite{NEOS1,NEOS2,NEOS3}.

\begin{table}[!t]
\caption{ Parameters of PSs Generation Units}
\centering
\begin{tabular}{ c c c c c c }
\hline
\hline
Ship		& 	a	& 		b 	& 		c &  SU/SD(\$) & Sailing(\$/h) \\
\hline
PS 1 	& 	74.33	 & 		15.4708		& 0.045923 &		45	&	250  \\
\hline
PS 2   & 	58.810 	& 		22.942				& 	0.00977 &		45	&	100   \\
\hline
\end{tabular}
\label{tgcps}
\end{table}

\subsection{6-Bus Power System Case Study}
\subsubsection{Case I}
\begin{figure}[!t]
\centering
\includegraphics[width=3.0in]{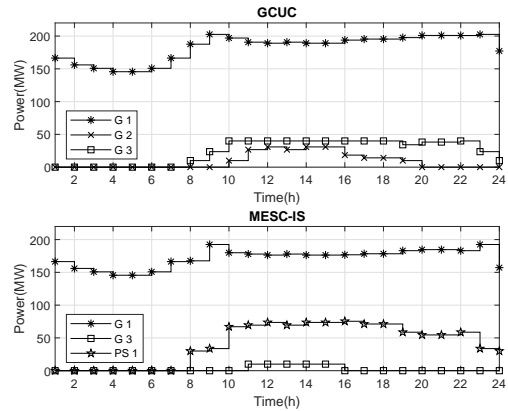} 
\caption{Generation units and PSs output in Case I (GCUC and MESC-IS).}
\label{figCI}
\end{figure}

\begin{figure}[!t]
\centering
\includegraphics[width=3.0in]{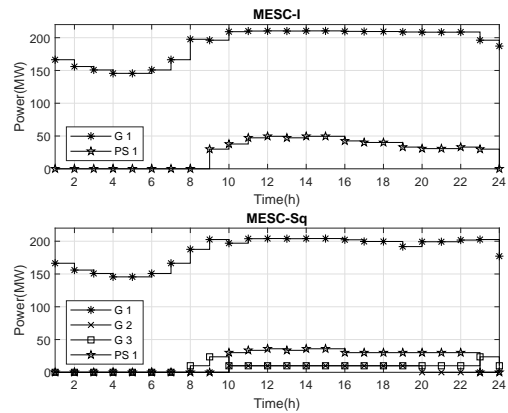} 
\caption{Generation units and PSs output in Case I (MESC-I and MESC-Sq).}
\label{fig2CI}
\end{figure}

First, the energy problem for the three islands was solved without considering the ships. Assuming deterministic load, the operational cost was \$ 87154.47, and load demand was fully met. Then, the MESC-IS model was solved with the PSs considered stationary, and the problem solved again. The new cost for the system was \$ 79133.8. That represented a significant cost reduction even when the ships were considered stationary during the planning horizon. Fig. \ref{figCI} and Fig. \ref{fig2CI} show the generation units and PSs output in both cases. PS 1 operated for 15 hours, from hour 8 to the end of the planning horizon, supplying electricity from the same bus at which it was initially located (Bus 3). The potential of further savings was examined by considering the other two solution approaches. Since we already had the unit status from the solution of the original model (GCUC), the sequential model(MESC-Sq) was considered first. MESC-Sq is the case when the ships are coordinated after solving the GCUC. For MESC-Sq, the total cost was \$ 85583.39. Since the coordination was made after fixing the generation unit status, all the generation units' status remained on. However, the reduction was due to the mobility of the ships. The ship PS 1, traveled from bus 2 to bus 4 in this case, departing at hour 6 and arriving at hour 10. Because we restricted the ships' ability to operate upon arrival, it started running at hour 10.
The optimal approach is when the problem is solved all at once. The fully integrated model (MESC-I) provided an optimal solution of \$ 78463.17. PS 1, in this case, made the same trip as it did in the solution of previous approach, but with a different departure time. The ship departed one hour earlier, to arrive at hour 9. A significant reduction was achieved with respect to the original solution; Table \ref{tCI} summarizes the costs. Another interpretation of the results can be stated as follows: Ships are still viable to improve the power system economical operations even if their operational costs (sailing, entering, and departure costs) are \$9,000 higher per day. However, a more extensive case would be more accurate to actual real-life implementation.
\begin{table}[!t]
\caption{ Ports Costs}
\centering
\begin{tabular}{c c c c}
\hline
 Power Ship	 	    & Waiting 		& Entering            & Berthing     \\
\hline
\hline
PS 1		    	    & 	55			& 		200	    &  	 235		\\
PS 2 	   	    	    & 	20 			& 		200    	     &  210		\\
\hline
\end{tabular}
\label{tpc}
\end{table}
\begin{table}[!t]
\caption{Total Costs Case I }
\centering
\begin{tabular}{c c c c }
\hline
GCUC(\$)	 &  MESC-IS(\$) & MESC-Sq(\$)	   & MESC-I(\$)        \\
\hline
\hline
87154.47   	    & 	79133.8	& 		85583.39 	&   78463.17  \\
\hline
\end{tabular}
\label{tCI}
\end{table}

\subsection{118-Bus Power System Case Study}
The 118-bus power system was used to simulate a larger case using the same ships from the previous case. Two different loading scenarios were considered, classified as  high and low loading levels. The high loading level had an hourly average load of 1366.83 MW, while the low loading had an average load of 1354.04 MW. However, the complete load data are avialable from \cite{data}.

\subsubsection{Case II}
\begin{figure}[!t]
\centering
\includegraphics[width=3.0in]{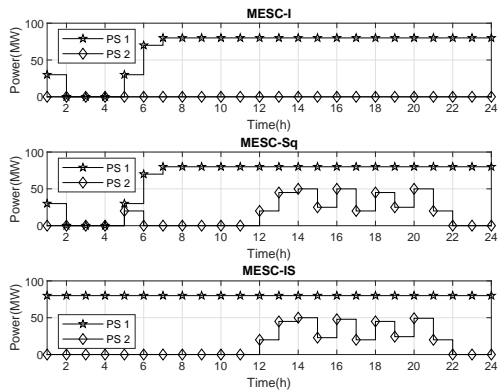} 
\caption{PSs generation output in case II.}
\label{figCII}
\end{figure}

In this case, the GCUC was solved for the system, providing a total cost of \$2052910.22 with 4.276 MW load shedding. In a similar manner to that considered in the previous case, three solution approaches were obtained. Table \ref{tCII} shows the total cost of the three approaches. The MESC-IS solution resulted in 3.86 MW shedding, a slight improvement compared to the original case. Even better, the MESC-I and MESC-Sq solutions resulted in no load shedding and offered significant cost reductions,of \$10696.6 and \$15906.12, respectively. Moreover, the mobility of the ships, in this case, contributed more in cost savings in comparison to Case I. As shown in Fig. \ref{figCII} PS 1 made one trip and generated 1570 MWh during the planning horizon; the ship traveled from its initial port to port 7. MESC-I suggested the deployment of PS 2. In fact, PS 2 made more trips than PS 1. MESC-I called for PS 2 to sail from bus 97 at the beginning of hour 3, and causing it to arrive at bus 10 at hour 5. It operated for only one hour and waited 4 hours before sailing to bus 7. Table\ref{tCIICIII} shows the complete coordination plan of the ships for the MESC-I and MESC-Sq approaches. 

\begin{table}[!t]
\caption{Total Costs Case II}
\centering
\begin{tabular}{c c c c }
\hline
GCUC(\$)	 &  MESC-IS(\$) & MESC-Sq(\$)	   & MESC-I(\$)      \\
\hline
\hline
2052910.22   	    & 	2043670.36	& 		2042213.61 	&   2037004.1 \\
\hline
\end{tabular}
\label{tCII}
\end{table}

\subsubsection{Case III}
\begin{figure}[!t]
\centering
\includegraphics[width=3.0in]{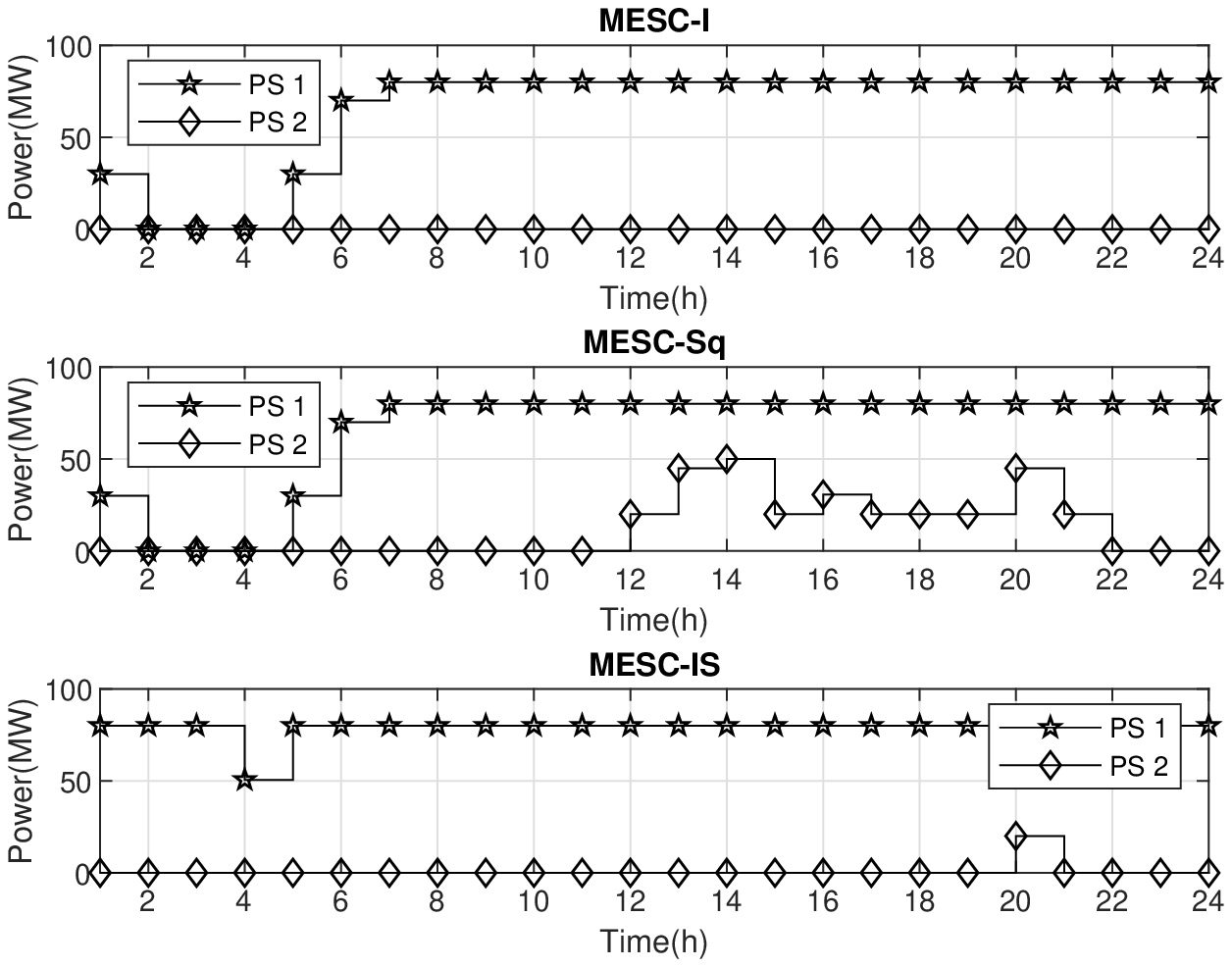} 
\caption{PSs generation output in case III.}
\label{figCIII}
\end{figure}

To further illustrate the potential of integrating the PSs, different load demand was considered in this case. The GCUC resulted in a total cost of \$2011241.41 with no load shedding. The Integrated solution resulted in a system savings of \$6971.68, while the sequential solution approach provided \$3964.8 in savings. Both ships sailed once during the planning horizon. As shown in Fig. \ref{figCIII}, the total energy produced by the ships during this time was 1860.7 MWh. The total operational costs of the ships were \$27605.086 and \$35328.263 for the sequential and the integrated approaches, respectively. In this case, the stationary solution approach resulted in a better solution than the sequential one. The MESC-IS total cost was \$2005790.52. Table \ref{tCIII} shows the costs of all solution approaches. Moreover, with relatively low demand and with no shedding penalties imposed, the ships were still able to show a potential for economic improvements.
\\
\\

\begin{table}[!t]
\caption{Total Costs Case III}
\centering
\begin{tabular}{c c c c }
\hline
GCUC(\$)	 &  MESC-IS(\$) & MESC-Sq(\$)	   & MESC-I(\$)          \\
\hline
\hline
2011241.41  	    & 2005790.52		& 		2007276.61 	&   2004269.73 \\
\hline
\end{tabular}
\label{tCIII}
\end{table}

\begin{table*}[!t]
\caption{ Results for Case II and Case III}
\centering
\begin{tabular}{c c c c c | c c c c c || c c c c c | c c c c c c c}
\hline
&\multicolumn{8}{c}{Case II}&&&\multicolumn{8}{c}{Case III}&\\
\hline
&\multicolumn{4}{c}{Integrated}&\multicolumn{4}{c}{Sequential}&&&\multicolumn{4}{c}{Integrated}&\multicolumn{4}{c}{Sequential}&\\
Time & \multicolumn{2}{c}{L}&\multicolumn{2}{c}{P } & \multicolumn{2}{c}{L}&\multicolumn{2}{c}{P } &&& \multicolumn{2}{c}{L}&\multicolumn{2}{c}{P } & \multicolumn{2}{c}{L}&\multicolumn{2}{c}{P} \\
&PS1&PS2&PS1&PS2&PS1&PS2&PS1&PS2&&&PS1&PS2&PS1&PS2&PS1&PS2&PS1&PS2&\\
\hline
\hline
1&87&97&30&0&87&97&30&0&&&87&97&30&0&87&97&30&0&\\

2&$|$&97&0&0&$|$ &97&0&0&&&$|$&97&0&0&$|$&97&0&0&\\

3&$|$&$|$&0&0&$|$&97&0&0&&&$|$&97&0&0&$|$&97&0&0&\\

4&$|$&$|$&0&0&$|$&97&0&0&&&$|$&97&0&0&$|$&97&0&0&\\

5&7&10&30&20&7&97&30&0&&&7&97&30&0&7&97&30&0&\\

6&7&10&70&0&7&97&70&0&&&7&97&70&0&7&97&70&0&\\

7&7&10&80&0&7&97&80&0&&&7&97&80&0&7&97&80&0&\\

8&7&10&80&0&7&97&80&0&&&7&97&80&0&7&97&80&0&\\

9&7&10&80&0&7&97&80&0&&&7&97&80&0&7&97&80&0&\\

10&7&$|$&80&0&7&97&80&0&&&7&$|$&80&0&7&97&80&0&\\

11&7&$|$&80&0&7&97&80&0&&&7&$|$&80&0&7&97&80&0&\\

12&7&7&80&20&7&97&80&0&&&7&7&80&20&7&97&80&0&\\

13&7&7&80&45&7&97&80&0&&&7&7&80&45&7&97&80&0&\\

14&7&7&80&50&7&97&80&0&&&7&7&80&50&7&97&80&0&\\

15&7&7&80&25&7&97&80&0&&&7&7&80&20&7&97&80&0&\\

16&7&7&80&50&7&97&80&0&&&7&7&80&30.67&7&97&80&0&\\

17&7&7&80&20&7&97&80&0&&&7&7&80&20&7&97&80&0&\\

18&7&7&80&45&7&97&80&0&&&7&7&80&20&7&97&80&0&\\

19&7&7&80&25&7&97&80&0&&&7&7&80&20&7&97&80&0&\\

20&7&7&80&50&7&97&80&0&&&7&7&80&45&7&97&80&0&\\

21&7&7&80&20&7&97&80&0&&&7&7&80&20&7&97&80&0&\\

22&7&7&80&0&7&97&80&0&&&7&7&80&0&7&97&80&0&\\

23&7&7&80&0&7&97&80&0&&&7&7&80&0&7&97&80&0&\\

24&7&7&80&0&7&97&80&0&&&7&7&80&0&7&97&80&0&\\
\hline
Total (MWh)&\multicolumn{4}{c}{1940}&\multicolumn{4}{c}{1570}&&&\multicolumn{4}{c}{1860.7}&\multicolumn{4}{c}{1570}&\\
\hline
PSs Total cost (\$)&\multicolumn{4}{c}{37850.536}&\multicolumn{4}{c}{27605.086}&&&\multicolumn{4}{c}{35328.263}&\multicolumn{4}{c}{27605.086}&\\
\hline
System Saving (\$)&\multicolumn{4}{c}{15906.12}&\multicolumn{4}{c}{10696.6}&&&\multicolumn{4}{c}{6971.68}&\multicolumn{4}{c}{3964.8}&\\
\hline
\end{tabular}
\label{tCIICIII}
\end{table*} 
\section{Computational Performance}

The generation scheduling problem has long been known as a computationally challenging to solve in real time. In fact, it is classified as a non-deterministic polynomial time hardness problem (NP-hard). Therefore, it is critical to test the computational performance of the problem whenever any changes are made. A significant computational burden was noticed with the integrated solution approach, MESC-I. A 2500 seconds running time was considered for all case studies. In Case II, the solver was able to achieve a global optimum solution for GCUC in 327.04 s and 1059585 iterations. However, in the integrated approach, MESC-I, the solver was unable to reach the global optimum and terminated after 2500 s and 2449504 iterations. In Case III, GCUC needed 410.6 s and 1007838 iterations, while MESC-IS and MESC-Sq required just 61.51 s and 24 s ,respectively, to reach a global optimum solution. Similar to the previous case, the MESC-I approach was unable to reach the global optimum within the given 2500 s, a period more than 104 times the time that required for sequential approach to achieve solution. Even though the results of the integrated approach are unsatisfying, they offer insight into the exponential growth in the complexity of the problem. Regardless, the approximate approaches performed computationally well, as shown in table \ref{tCOMII} and table \ref{tCOMIII}, leaving a time-accuracy trade-offs potential between the different solution approaches. 

It is important to mention that the MESC-Sq solution time in table \ref{tCOMIII} and table \ref{tCOMIII} only include the maritime-based problem. The first sub-problem must be solved, even without the existence of PSs. The total time or wall time can be obtained with a high degree of accuracy by adding the times together. Since we assume the availability of the energy scheduling solution, the two times are provided independently.
\\
\\

\begin{table}[!t]
\caption{Computational Complexity Case II }
\centering
\begin{tabular}{c c c c c}
\hline
&	   Objective (\$) & CPLEX (Sec)	   & Iteration    &  R-Gap \%  \\
\hline
\hline
GCUC 	    & 	2052910.22 &327.04 &1059585				&   0  \\
MESC-IS	  	    & 	2043670.36&524.52& 745422			 	&   0  \\
MESC-Sq  	    & 	2042213.61 &50.82&89206		 	&   0 \\
MESC-I   	  	    & 	2037004.1&2500&	2449504 	    	& 0.000163    \\
\hline
\end{tabular}
\label{tCOMII}
\end{table}
 
\begin{table}[!t]
\caption{Computational Complexity Case III }
\centering
\begin{tabular}{c c c c c}
\hline
&	   Objective (\$) & CPLEX (Sec)	   & Iteration    &  R-Gap \%  \\
\hline

\hline
GCUC 	    & 	2011241.41 &410.62&1007838				&   0  \\
MESC-IS  	    & 	2005790.52&61.51& 118964			 	&   0  \\
MESC-Sq 	    & 	2007276.61 &24.12&27709		 	&   0 \\
MESC-I 	  	    & 	2004269.73&2500&	1814246 	& 0.000076    \\
\hline
\end{tabular}
\label{tCOMIII}
\end{table}

\section*{Conclusion}
In this paper, a model was presented to incorporate maritime energy transportation models into power system operation. The model effectively coordinated power-generating ships with the electric grid and demonstrated an economic saving potential under different scenarios. 
On the other hand, computational complexity of the problem was a major concern. As discussed, the complexity of the problem was exponentially increasing with the problem size. However, since this is one of the first papers to consider such incorporation, there is a potential to introduce valid inequalities and cuts to improve the computational performance of the model.

\ifCLASSOPTIONcaptionsoff
  \newpage
\fi
\bibliographystyle{IEEEtran}
\bibliography{ref}
\end{document}